\documentclass[11pt]{amsart} 

\usepackage[utf8]{inputenc} 
\usepackage[T1]{fontenc}
\usepackage{lmodern}

\usepackage{amssymb}
\usepackage{amsmath}
\usepackage{mathtools}
\usepackage[mathscr]{euscript}
\usepackage{multicol}

\usepackage{geometry} 
\usepackage{graphicx} 

\usepackage{booktabs} 
\usepackage{array} 
\usepackage{paralist} 
\usepackage{verbatim} 
\usepackage{subfig} 
\usepackage{xcolor}

\usepackage{fancyhdr} 
\newtheorem{lemma}{Lemma}

\newtheorem{cor}[lemma]{Corollary}
\newtheorem{thm}[lemma]{Theorem}
\newtheorem{defn}[lemma]{Definition}
\newtheorem{conj}[lemma]{Conjecture}

\newtheorem{rmk}[lemma]{Remark}
\newtheorem{example}[lemma]{Example}

\newcommand{\Z}{\mathbb Z}
\newcommand{\R}{\mathbb R}
\newcommand{\C}{\mathbb C}

\newcommand{\po}{\mathbb P^1}
\newcommand{\ql}{\mathbb Q_\ell}
\newcommand{\zl}{\mathbb Z_\ell}

\newcommand{\Q}{\mathbb Q}
\newcommand{\qb}{\bar{\mathbb Q}}
\newcommand{\F}{\mathbb F}

\newcommand{\fl}{\mathbb F_\ell}

\newcommand{\kb}{\bar k}

\newcommand{\kv}{k_v}

\newcommand{\zv}{\mathbb Z_v}

\newcommand{\ok}{\mathcal O_k}
\newcommand{\ov}{\mathcal O_v}
\newcommand{\pp}{\mathfrak p}
\newcommand{\ml}{\mathfrak L}
\newcommand{\PP}{\mathfrak P}

\newcommand{\dc}{\mathfrak D_C}
\newcommand{\s}{\mathcal S}

\newcommand{\ang}[1]{\langle #1 \rangle}

    \DeclareFontFamily{U}{wncy}{}
    \DeclareFontShape{U}{wncy}{m}{n}{<->wncyr10}{}
    \DeclareSymbolFont{mcy}{U}{wncy}{m}{n}
    \DeclareMathSymbol{\Sha}{\mathord}{mcy}{"58}

\title{Hasse Principle Violations in Twist Families of Superelliptic Curves}
\author{Lori D. Watson}

\begin{document}
\maketitle

\begin{abstract}
Conditionally on the $abc$ conjecture, we generalize previous work of Clark and the author to show that a superelliptic curve $C: y^n = f(x)$ of sufficiently high genus has infinitely many twists violating the Hasse Principle if and only if $f(x)$ has no $\Q$-rational roots. We also show unconditionally that a curve defined by $C: y^{pN}=f(x)$ has infinitely many twists violating the Hasse Principle over any number field $k$ such that $k$ contains the $p$th roots of unity and $f(x)$ has no $k$-rational roots.
\end{abstract}

\section{Introduction}
Let $k$ be a number field and $C/_k$ superelliptic curve\footnote{All curves are assumed to be \textit{nice}: that is smooth, projective, and geometrically integral.} (see Definition \ref{supdef}).
The main result of this work is that infinitely many twists\footnote{Given a curve $C$ defined over $k$, a twist of $C$ is a curve $C'$ also defined over $k$ such that $C$ and $C'$ become isomorphic over $\kb$.} of a superelliptic curve will violate the Hasse Principle unless the existence of $\Q$-rational roots of $f$ or the degree of $f$ guarantees a $\Q$-rational point on every twist. 
\begin{defn}\label{supdef}
Let $n \ge 2$ be an integer and let $k$ be a  number field. A superelliptic curve $C/_k$ is a nice curve admitting an affine model of the form $y^n = f(x)$, where $f(x) \in k[x]$ factors over $\kb$ as 
\\ \centerline{$\displaystyle f(x) = A \prod_{j=1}^r(x-\alpha_j)^{n_j},$}
where $A \in \kb^*$,  the $\alpha_j$ are distinct, $1 \le n_j < n$ for each $j$, and $\gcd(n, n_1, \ldots, n_r) = 1$. 
\end{defn}
\begin{rmk}
As shown in \cite[Lemma 1]{KOO}, the condition $\gcd(n, n_1, \ldots, n_r) = 1$ guarantees that $C$ is geometrically irreducible.
\end{rmk} 

\begin{rmk}
When $n=2$, the curve is hyperelliptic.
\end{rmk}

 Over an algebraic closure $\kb$ of $k$, such a curve is equipped with an automorphism $\tau$ of order $n$ defined by $\tau(x,y) = (x, \zeta y)$, where $\zeta$ is a primitive $n$th root of unity. One especially interesting feature of these curves is that their automorphism groups give rise to families of twists.  In this work, we consider twists $C_d: dy^n = f(x)$, where $d$ is an $n$th power free integer. Though $C$ and $C_d$ define the same geometric object, their arithmetic can be very different. 
\\
\begin{example}
The elliptic curve $E: y^2 = x^3+7x$ has finitely many $\Q$-rational points while some of its quadratic twists have infinitely many, including $E_5: 5y^2=x^3+7x$ (see \cite[P. 351]{AEC}).
\end{example}
In Theorem \ref{main} we show that if one assumes the $abc$ conjecture, then a family of twists $C_d$ of a superelliptic curve $C: y^n = f(x)$ exhibits one of two behaviors: either there are no Hasse Principle violations within the family of twists, or there are many. Follwing Mazur and Rubin (\cite{MR}), for a given family of twists we say that ``many'' curves in the family exhibit property $P$ if  the number of twists $C_d$ such that $|d| \le X$ and $C_d$ satisfies property $P$ is $\gg X/(\log X)^\gamma$ for some $\gamma \in \R$. In order for the family of twists to have infinitely many Hasse Principle violations (or indeed any), it is necessary that the map $C \rightarrow C/\langle \tau \rangle \cong \po$ has no $\Q$-rational branch points; assuming the $abc$ conjecture, we show that this condition is also sufficient. In particular, we show

\begin{thm}
\label{main} Assume the $abc$ conjecture. Let $C: y^n = f(x)$ where $n \ge 2$ is an integer and $f(x) \in \Z[x]$ has distinct roots in $\bar \Q$ and $\deg(f) \ge \frac{4n-2}{n-1}$. Then the following are equivalent:
\\ [0.1\baselineskip] \indent (i) The automorphism $\tau$ has no $\Q$-rational fixed points.
\\ [0.1\baselineskip] \indent  (ii) There exists $0 \le \gamma < 1$ such that as $X \rightarrow \infty$, the number of $n$th power free integers $d$ \\ \indent with $|d| \le X$ such that $C_d$ violates the Hasse Principle is $\displaystyle \gg_C \frac{X}{\log (X)^\gamma}$.
\\ [0.1\baselineskip] \indent  (iii) Some degree $n$ twist $C_d$ violates the Hasse Principle.
\end{thm}
That (ii) implies (iii) is immediate. For (iii) implies (i), a $\Q$-rational branch point of $\ang \tau$ will remain rational on every twist. Thus it remains to prove (i) implies (ii), which we do in Section \ref{smain}. For this, we will require two results. The first is due to Granville (\cite{G}) and it provides an upper bound on the number of twists that have so-called ``nontrivial'' $\Q$-rational points. By nontrivial, Granville means a point that is \emph{not} a fixed point of the `superelliptic automorphism. In the absence of any $\Q$-rational branch points, Granville's results provides an upper bound on the number of twists which have any $\Q$-rational points at all. The second result is a strengthening of \cite[Thm. 3]{CW}. It provides a lower bound on the number of twists of a curve that have points everywhere locally (and in fact shows that there are many such twists).

This result generalizes Theorem 2 of \cite{CW}, which proves the result when $C$ is hyperelliptic. In \cite{KL}, K\"{o}nig and Legrand recover the hyperelliptic result. They then extend the result to Galois covers of $\po/_\Q$, that is, when the Galois group of $\Q(C)/\Q(\po)$ contains in its center $\Z/n\Z$. We note that this result of  K\"{o}nig and Legrand is separate from Theorem \ref{main} except in the hyperelliptic case: while the Galois group of $\qb(C)/\qb(\po)$  does contain $\ang \tau \cong \Z/n\Z$ for a superelliptic curve $C$, such a curve gives a Galois cover of $\po/_\Q$ only when $n = 2$. It gives a Galois cover of $\po_k$ if and only if $k$ contains the $n$th roots of unity. In that sense, we may think of superelliptic curves as geometrically Galois covers of $\po$.
\par In Section \ref{smain} we  give the proof of Theorem \ref{main}. In Section \ref{sab} we provide asymptotic bounds for the number of Hasse Principle violations. In Section \ref{uncon}, we give unconditional results on Hasse Principle violations over a number field containing $p$th roots of unity (where $p$ is prime) for curves of the form $y^{pN}=f(x)$.

\section{Proof of the Main Result}
\label{smain}
In this section  we prove the main result. 
\par We begin by providing criteria for the existence of Hasse Principle violations within a family of twists. A point $P$ fixed under the action of the automorphism group $\ang \tau$ is either a point lying above $\infty$ or a point with $y$-coordinate $0$. The quotient map $C \rightarrow C/\ang \tau \stackrel{\sim}{\rightarrow} \po$ has a $\bar \Q$ branch point above $\infty$ if and only if $n \nmid \deg(f)$ with the ramification index of the point at infinity being $e = n/\gcd(n,\deg(f))$ (\cite{KOO}). 
A point $P$ with $y$-coordinate $0$ necessarily has as its $x$-coordinate a root of the defining polynomial. 

\subsection{Bound on twists with nontrivial global points}
The ``global'' step towards proving the main theorem is a result due to Granville. This results provides an upper bound on the number of twists with nontrivial $\Q$-rational points. It relies on the $abc$ theorem to bound the height of the $x$-coordinate of a point that can appear on a twist $C_d$.
\subsection{The $abc$ conjecture}
One of several equivalent statements of the $abc$ conjecture over $\Q$ is as follows (\cite[Conj. 2.4]{V}):
\begin{conj}[$abc$ Conjecture, Masser-Oesterl{\'e}]
For each $\epsilon >0$, there exists a constant $K$ depending only on $\epsilon$ such that for all triples $(a,b,c)$ of relatively prime  integers with $a+b = c$,
\\ \centerline{$\displaystyle \max\left(|a|, |b|, |c|\right)\le K( \prod_{\substack{p \mbox{\tiny{ prime}} \\ p|abc}} p)^{1+\epsilon}$,}
\end{conj}
The $abc$ conjecture has numerous implications; most relevant to this work is that the $abc$ conjecture, if true, allows one to count the number of $n$th power-free values achieved by a polynomial with integer coefficients. 
\begin{thm}[Granville, {\cite[\S 11]{G}}]\label{gmain}
Assume that the $abc$-conjecture is true. Fix $\epsilon > 0$. Let $f(x) \in \Z[x]$ have no repeated roots in $\qb$. Let $d \in \Z$ be an $n$th power free integer. Any rational point on $C_d: dy^n = f(x)$ with $x$-coordinate $r/s$ where $\gcd(r,s) = 1$ satisfies
\\ \centerline{$\displaystyle |r|, |s| \ll_{f,\epsilon} |d|^{-\left(nk+i-1 - \frac{\gcd(n,i)+1}{n-1}\right)+\epsilon}$,}
where $\deg(f) = nk+i$, $1 \le i \le n$. \end{thm}

As a consequence, we conclude that as $X \rightarrow \infty$, the number of $n$th power free integers $d$ with $|d| \le X$ for which $C_d$ has a nontrivial rational point is $\displaystyle \ll_f X^{2/\left(nk+i-1 - \frac{\gcd(n,i)+1}{n-1}\right)+\epsilon}$

If $\frac{4n-2}{n-1}< \deg(f)$, then for sufficiently small $\epsilon >0$, the number of $n$th power free $d \in \Z$ with $|d| \le X$ such that $C_d$ has a nontrivial rational point is $X^a$ for some $a < 1$.  
\subsection{Bound on twists with points everywhere locally}
With an eye towards eventually proving the main theorem for arbitrary number fields, in this section we prove that if a superelliptic curve $C$ has points everywhere locally, then many of its twists has points everywhere locally as well.
\\ 
\begin{lemma} \label{locbnd} Let $C: y^n = f(x)$ be a superelliptic curve of genus $g \ge 1$ defined over a number field $k$. Suppose $C$ has points everywhere locally. For an element $d \in \ok$, let $N(d) = \vert \ok/(d) \vert$. As $X \rightarrow \infty$, the number of non-associate elements $\pi \in \ok$ such that $N(\pi) \le X$ and $C_\pi$ has points everywhere locally is  $\gg_{C,k} X/\log X$.
\end{lemma}

\textbf{Proof}:  Let $(\pi) = \pp \triangleleft \ok$ denote a principal maximal ideal of $\ok$ with $\pi$ satisfying $\sigma(\pi) > 0$ for every real embedding $\sigma: k \hookrightarrow \R$. (These ideals are precisely those which split in the narrow Hilbert Class Field of $k$, so by the finiteness of the class number they represent a positive density set of the prime ideals of $\ok$, (\cite[\S 3]{J})). For any place $v $ of $k$, if $\pi \in \kv^{*n}$, then $C_\pi \cong C$ over $\kv$. Thus, if $k_v$ is any archimedean completion of $k$ then $C_\pi \cong C$ over $\kv$and $C_\pi(\kv) \neq \emptyset$. From now on, $v$ will denote a finite place corresponding to a prime ideal $\mathfrak L_v$, $\ov$ will denote the ring of integers of $\kv$, and $\fl$ the finite field with $\ell$ elements will denote the residue field of $\kv$.
\par Let $M_1 \in \Z^+$ be such that $C$ extends to a smooth proper relative curve over $\ov$ for every $v$ such that $\operatorname{char}(\mathfrak L_v) > M_1$. Such an $M_1$ exists for any nice curve $C/_{k}$ by the openness of the smooth locus. 
\par Suppose $\operatorname{char}(\mathfrak L_v)  > M \coloneqq \max\{M_1, 4g^2-1, n\}$, $\ml_v \neq \pp$, and $\pi \not \in \kv^{*n}$. Then the minimal regular model $C/_{\ov}$ is smooth. Fix an extension $\kv\left(\sqrt[n]{\pi}\right)/\kv$ which contains an $n$th root of $\pi$. Over $\kv\left(\sqrt[n]{\pi}\right)$, $C \cong C_\pi$. Since $\kv\left(\sqrt[n]{\pi}\right)/\kv$ is unramified for any choice of an $n$th root of $\pi$, and formation of the minimal regular model commutes with \'etale base change (\cite[Prop. 10.1.17]{LIU}), it follows that the minimal regular model $(C_\pi)_{/\ov}$ is smooth. By the Riemann hypothesis for curves over a finite field (\cite[Thm. V.2.2]{AEC}), since $\ell \geq 4g^2$, we have $C_\pi(\F_{\ell}) \neq \emptyset$, so by Hensel's Lemma (\cite{Cass})  we have $C_\pi(\kv) \neq \emptyset$.
\par Suppose now $\operatorname{char}(\mathfrak L_v)  \le  M$ and $\ml_v \neq \pp$. If the $v$-adic valuation of $n$ is $\theta$, then we require $\pp$ to split completely in  $\tilde k = k(\zeta_{\ell^{2\theta+1}})$, the compositum of $k$ and $\Q(\zeta_{\ell^{2\theta+1}})$. In that case we have $\pi \equiv 1 \pmod{\ell^{2\theta+1}}$; letting $g(x) = x^n-\pi$, we then have $|g(1)|_v  \le \ell^{-(2\theta+1)} < \ell^{-2\theta} = |g'(1)|_v^2$, so by Hensel's lemma, there exists a unique $u \in \zv$ such that $g(u) = 0$. Thus $\pi$ is an $n$th root in $\kv$,  so  $C/_{\kv} \cong C_\pi/_{\kv}$, and $C_\pi(\kv) \neq \emptyset$. 
\par Now suppose $\ml_v = \pp$. Let $P \in C(\bar k)$ be a superelliptic branch point. We assume $\pp$ splits completely in $K = k(P)$. Then, if $\PP$ is a prime of $K$ lying above $\pp$, since $\pp$ splits completely, the completion $K_{\PP}$ has $[K_{\PP}: k_\pp] = 1$, thus $K$ embeds into $k_\pp$, and $P \in  C_\pi(k_\pp)$, so $  C_\pi(k_\pp) \neq \emptyset$. We have imposed finitely many conditions on $\pp$, each requiring that $\pp$ splits completely in a certain number field. Letting $L$ be the Galois closure of the compositum of these finitely many number fields, we have that $C_\pi$ has points everywhere locally whenever $(\pi) = \pp$ splits completely in $L$. By the Chebotarev density theorem (\cite[Appendix]{LS}), this set of primes (which we will denote by $S$ and use in the next theorem) has positive density in the set of $\ok$ primes. By Landau's Prime Ideal Theorem, the number of prime ideals $\pp$ of $\ok$ with $N\pp \le X$ is asymptotic to $X/\log X$.  \hfill $\square$
\\$\;$
\par In the case $k = \Q$, having produced a positive density set of primes whose twists have points everywhere locally, we can construct a larger set of $n$th power free integers $d \in \Z$ such that the twists  $C_d$ have points everywhere locally.
\begin{thm}\label{locs} Let $C: y^n = f(x)$ be a superelliptic curve of genus $g \ge 1$ defined over $\Q$. Suppose $C$ has points everywhere locally. The number of $n$th power free integers $d$ with $|d| \le X$ such that $C_d$ has points everywhere locally is $\gg_C X/\log(X)^\gamma$ for some $0 \le \gamma < 1$.
\end{thm}

\textbf{Proof}: Denote the set of primes $p$ constructed in Lemma \ref{locbnd}  by $S$. By construction this set  has  density $0 < \delta(S) < 1$. Consider the set $\mathcal D$ consisting of $n$th power free integers $d$, all of whose prime divisors are in $S$. If $n$ is even  we further require that all elements of $\mathcal D$ be positive. We will show that for each $d \in \mathcal D$, the twist $C_d$ has points everywhere locally. 
\par If $n$ is even  then as $d > 0$ and $C(\R) \neq \emptyset$, we have $C_d(\R) \neq \emptyset$. If $n$ is odd or if $f(x)$ has a real root then $C_d(\R)  \neq \emptyset$ for every nonzero integer $d$.
\par Next, let $M$ be as in Lemma \ref{locbnd}, and let $\ell > M$, $\ell \nmid d$. If $d \in \Q_{\ell}^{* n}$, then $C_d$ is isomorphic to $C$ over $\Q_\ell$, thus $C_d(\Q_\ell) \neq \emptyset$. So assume $d \notin \Q_{\ell}^{* n}$. Then $\Q_\ell(\sqrt[n]{d})/\Q_\ell$ is unramified, and as before we conclude that the minimal regular model $(C_d)_{/\Z_{\ell}}$ is smooth. Then by the Riemann Hypothesis for curves over a finite field, since $\ell \ge 4g^2$, $C_d(\mathbb F_\ell) \neq \emptyset$, so by Hensel's Lemma, we again have $C_d(\Q_\ell) \neq \emptyset$.
\par For $\ell < M$ and $\ell \nmid d$, as $d = p_1^{m_1}\cdots p_r^{m_r}$ with $p_j \in \Q_{\ell}^{* n}$ for all $1\le j\le r$, $d \in \Q_{\ell}^{\times n}$, so $C_d$ is isomorphic to $C$ over $\Q_\ell$, thus $_d(\Q_\ell) \neq \emptyset$.
\par Finally, consider $C_d(\Q_p)$, where $p \vert d$. Then $p \in S$. By construction, for each $p \in S$, $p$ splits completely in $K = \Q(P)$ where $P$ is a fixed point of$\tau$. As before, if $\PP$ is a prime of $K$ lying above $p$, since $p$ splits completely, the completion $K_{\PP}$ has $[K_{\PP}: \Q_p] = 1$, thus $K$ embeds into $\Q_p$. Since $P$ is a fixed point of $\tau$, it is of the form $P = (\alpha, 0)$, where $\alpha$ is a root of the defining polynomial $f(x)$. Thus $K = \Q(\alpha) \subset \Q_p$. Then $(\alpha, 0)$ is a $\Q_p$-rational point of every degree $n$ twist of $C$, so in particular $C_d(\Q_p) \neq \emptyset$. Thus $C_d$. has points everywhere locally for each $d \in \mathcal D$.
\par Let $\gamma = 1-\delta(S)$. By \cite[Thm 2.4]{SD}, we have that the number of $d \in \mathcal D$  with $|d| \le X$ such that $C_d$ has points everywhere locally is $\gg X/\log (X)^\gamma$.  \hfill $\square$
\subsection{Proof of the Main Theorem}
$\;$\\
We now complete the proof of the main theorem.
$\;$ \\
\par \textbf{Proof of \ref{main}}: Let $C: y^n = f(x)$ be a superelliptic curve, where $f(x)$ has coefficients in $\Z$, distinct roots in $\qb$, no roots in $\Q$, and  $\deg(f) > {(4n-2)/(n-1})$. If $C$ does not have points everywhere locally, write $f(1) = d_1d_2^n$ where $d_1, d_2 \in \Z$ and $d_1$ is $n$th power free. Then the curve $C' \coloneqq C_{d_1}$ (with model $d_1y^n = f(x)$) has the $\Q$-rational point $(1, d_2)$ and hence has points everywhere locally. Applying Theorem \ref{locs} to $C'$, we have that the number of $n$th power free integers $d$ with $|d| \le X$ such that the twist $C'_{d'}$ of $C'$ with points everywhere locally is $\gg X/\log (X)^\gamma$. As the twist $C'_{d'}$ is a twist of the original curve $C$, we have that the number of $n$th power free integers $d$ with $|d| \le X$ such that $C_d$ has points everywhere locally is still $\gg_C X/\log (X)^\gamma$. 
\par By \ref{gmain}, we have that the number of $n$th power free $d \in \Z$ with $|d| \le X$ such that $C_d(\Q) \neq \emptyset$ is $\ll X^{2/3}$. Thus the number of $n$th power free $d$ with $|d| \le X$ such that $C_d$ violates the Hasse Principle is $\gg_C  X/\log (X)^\gamma$.   \hfill $\square$

\section{Asymptotic Bounds }\label{sab}
Assuming the $abc$ conjecture, Theorem \ref{main} shows that the existence of any (equivalently infinitely many) Hasse Principle violations within a family of twists of a curve $C: y^n = f(x)$ depends on the degree of $f$ and whether $f$ possesses any $\Q$-rational roots. By more closely examining the local behavior of $f$, we can provide more precise bounds on the number of twists violating the Hasse Principle.
\par Let $\dc \coloneqq \{n\mbox{th power free } d\in \Z : C_d \mbox{ has points everywhere locally}\}$.
\\ For $X \ge 1$, put
$\dc(X) = \#\dc\cap[-X, X]$. We saw in Thm \ref{locs}, that $\dc(X) \gg \frac{X}{\log (X)^\gamma}$ for some $\gamma \in (0,1)$. We will soon see that $\dc(X)$ (and hence the number of Hasse Principle violations for $C$ with no $\Q$-rational branch points) depends on the density (within the set of all primes) of the set $\s= \{\ell \mbox{ prime} : f(x) \mbox{ has a root modulo } \ell\}$. For the hyperelliptic case $n=2$, we can provide an unconditional upper bound on the number of twists having points everywhere locally in terms of the density $\delta$ of $\s$. For $n \ge 2$, when the density of $\s$ equals 1, we will show that, conditional on the $abc$ conjecture, a positive density set of twists have points everywhere locally. Before proceeding we introduce some new terminology:

\begin{defn}
We say a polynomial $f \in \Z[x]$ is \emph{weakly intersective} if $\delta = 1$.
\end{defn}

\begin{thm}\label{dc1}
Let $C: y^n = f(x)$ be a superelliptic curve with $f(x) \in \Z[x]$ squarefree and weakly intersective. Then $\dc(X) \gg_C X$.
\end{thm}

As an immediate consequence of Theorem \ref{dc1} we have the following corollary:

\begin{cor}
Let $C:y^n = f(x)$ be a superelliptic curve with $f(x) \in \Z[x]$ squarefree. If $C$ has no $\Q$-rational superelliptic branch points and $\deg(f) > 5$, then conditionally on the $abc$ conjecture, as $X \rightarrow \infty$, the number of degree $n$ twists of $C/_{\Q}$ that violate the Hasse Principle is $\gg_C X$.
\end{cor}
\begin{rmk}
Before proving Theorem \ref{dc1}, we will first show that if $f(x)$ is weakly intersective, then the set of primes $\ell$ for which $f(x)$ does not have a root $\pmod \ell$ is finite. (This argument appears in \cite{CW}). Let $f(x) = \sum_{j=1}^n a_jx^j \in \Z[x]$ have degree $n \ge 2$ and let $\Delta$ be the discriminant of $f$. Suppose $f(x)$ has distinct roots in $\qb$ and let $G$ denote the Galois group of $f$. 
\par For each prime $\ell \nmid a_n\Delta$, the partition of $n$ given by the cycle type of a Frobenius element $\sigma_\ell$ at $\ell$ coincides with a partition of $n$ given by the degrees of the irreducible factors of the image of $f$ in $(\Z/\ell \Z)[x]$. Since $f(x)$ is weakly intersective, by the Frobenius Density Theorem (\cite[\S 3]{LS}), every $\sigma \in G$ has a fixed point, and thus $f$ has a root $\pmod \ell$ for every $\ell \nmid a_n\Delta$. By Hensel's Lemma, $f(x)$ has a root in $\Z_\ell$ for all but finitely many $\ell$.
\end{rmk}
\noindent \textbf{Proof of \ref{dc1}}: By the remark above, if $f$ is weakly intersective, then $f$ has a root modulo $\ell$ for all but finitely many $\ell$ and hence has a root in $\zl$ for all but finitely many primes $\ell$. Therefore, the set $\mathcal P$ of primes $\ell$ such that $C(\ql) = \emptyset$ is finite. For each $\ell \in \mathcal P$, we have $C_d(\ql) \neq \emptyset$ for any $d$ lying in the same $\ql$-adic $n$th power class as $f(1)$. The set of integers lying in any given $\ql$-adic $n$th power class is a nonempty union of congruence classes modulo $\ell^{2v_\ell(n)}$ if $\ell $ is odd and modulo $2^{4v_\ell(n)}$ if $\ell = 2$. By the Chinese Remainder Theorem there are $a, N \in \Z^+$ such that if $d \equiv a \pmod N$, then $C_d(\ql) \neq \emptyset$ for all primes $\ell$. A result of Prachar (\cite{PRA}) guarantees that  there is a positive density set of $d \equiv a \pmod N$ which are squarefree (and thus $n$th power free), so long as $a \in \left(\Z/N\Z\right)^*$.  If $f$ has a real root, then $C_d(\R) \neq \emptyset$ for all $n$th power free $d \in Z$. Otherwise, $C_d(\R) \neq \emptyset \iff df(1) >0$. In either case, $\dc(X) \gg_f  X$. (The implied constant depends on the discriminant of $f$.) \hfill $\square$

\section{Unconditional Results}\label{uncon}
In this section, we provide an unconditional result on Hasse Principle violations in families of twists of certain superelliptic curves which map to curves of sufficiently large genus and few points. We prove a theorem, first stated in (\cite[Thm. 2]{TAHP1}), to produce many twists which fail to have $k$-rational points and thus, combined with Lemma \ref{locbnd} yield many Hasse Principle violations within some families of twists.
\begin{thm}\label{TAHP}
Let $k$ be a number field containing the $p$th roots of unity where $p$ is prime. Let $C/k$ be a smooth, projective, geometrically integeral curve, and let $\psi: C \longrightarrow C$ be a $k$-rational  automorphism of order $p$. Assume the following hold:
\par (i) $\{P \in C(k) : \psi(P) = P\} = \emptyset$.
\par (ii) $\{P \in C(\bar k) : \psi(P) = P\} \neq \emptyset$.
\par (iii) For some extension $L = k\left(d^{1/p}\right)$, the twist $C_d$ has points everywhere locally.
\par (iv) The set $(C/\ang \psi)(k)$ is finite.
\\ Then for all but finitely many $d$, the twisted curve $C_d$ has no $k$-rational points.
\end{thm}

\textbf{Proof}: Let $L/k$ be a cyclic degree $p$ extension. As $k$ contains the $p$th roots of unity, by Kummer Theory, $L = k\left(d^{1/p}\right)$, for some $p$th power free $d \in k$. Let $Y \coloneqq Y_L$ be the twist of $C$ by $\psi$ with respect to $L/k$. Then $\psi$ defines a $k$-rational automorphism on $Y$. Let $\theta: G_k \longrightarrow Aut(C)$ be the 1-cocycle corresponding to the twist $Y$. For a generator $\sigma$ of the Galois group of $L/k$, we have $\theta(\sigma) = \psi^j$ for some $j = j_\sigma \in \left(\Z/p\Z\right)^*$. We have that $Y/_L \cong C/_L$, and $\sigma$ acts on $Y(L)$ by $\sigma^*(P) = \psi^j\sigma(P)$. Thus, for all $P \in Y(A)$ (for any $k$-algebra $A$) we have
\\ \centerline{$\sigma^*\psi(\sigma^*)^{-1} = (\psi^j\sigma)\psi(\psi^j\sigma)^{-1} = (\psi^j\sigma)\psi(\sigma^{-1}\psi^{-j}) = \psi^j(\sigma\psi\sigma^{-1})\psi^{-j} = \psi^j\psi\psi^{-j} = \psi$}.
$\;$
We have natural maps
\\ \centerline{$\kappa: Y_L(k) \hookrightarrow C(L)$}
and
\\ \centerline{$\lambda: C(L) \rightarrow (C/\ang \psi)(L)$}
so that
\\ \centerline{$S_L \coloneqq (\lambda \circ \kappa)(Y_L(k)) \subseteq (C/\ang\psi)(k)$ and $(C/\ang\psi)(k) = \psi(C(k)) \cup  \bigcup_{L = k(d^{1/p})} S_L $.}
For distinct degree $p$ extensions of $k$, $L_1$ and $L_2$, $P \in S_{L_1} \cap  S_{L_2} \implies P \in C(L_1) \cap C(L_2) = C(k)$ (since $[L_i: k]$ is prime and the extensions are distinct, $L_1 \cap L_2 = k$). If $P \in S_L \cap C(k)$, then $P = \lambda(\kappa(Q))$ for some $Q \in Y_L(k)$. As $ \lambda(\kappa(Q)) \in (C/\ang \psi)(L)$, $P$ is a fixed point of $\psi$, but by hypothesis, no such points are $k$-rational. Thus $S_{L_1} \cap  S_{L_2} = \emptyset$ and $(C/ \ang \psi)(k)$ is a disjoint union of the $S_L$. Since $(C/ \ang \psi)(k)$ is finite, we conclude that there are only finitely many twists $Y_L$ which have $k$-rational points. \hfill $\square$
\begin{cor}\label{unsup}
Let $C:y^n = f(x)$ be a superelliptic curve defined over a number field $k$ with no $k$-rational superelliptic fixed points. Suppose $n = pN$ with $1< N < n$ and $p$ prime, and that $k$ contains the $p$th roots of unity. Let $N(d) \coloneqq |\ok/(d)|$ denote the norm of $d$. Suppose  that $Aut(C/{\kb}) \cong \mu_n$ and that the curve $C^{(N)}: y^N = f(x)$ has finitely many $k$-rational points. Then as $X \rightarrow \infty$ the number of $N$th power free $d \in \ok$ such that $C_d$ violates the Hasse Principle is $\gg_C \frac{X}{\log(X)}$. In particular, if the genus $g(C^{(N)}) \ge 2$ or if $C^{(N)}$ is an elliptic curve with finite Mordell-Weil rank over $k$, then $C$ has many twists violating the Hasse Principle.
\end{cor}

\textbf{Proof}: We assume without loss of generality that  $C$ has points everywhere locally over $k$. From Lemma \ref{locbnd}, there is a set of non-associate, totally positive prime elements $\pi \in \ok$ such that $C_\pi$ has points everywhere locally. We consider now the twists $C_d$ where $d = \pi^p$. As before, there is $M \in \Z^+$ (depending on $n$ and the genus and discriminant of $C$) such that for every prime ideal $\mathfrak L_v \neq (\pi)$, with $N\mathfrak L_v > M$, $C_d(k_v) \neq \emptyset$ (where $v$ is the finite place correspdoning to $\mathfrak L_v$, and $k_v$ the completion of $k$ at $v$). For $\mathfrak L_v \neq (\pi)$ with $N\mathfrak L_v \le M$, $\pi$, and hence $d$ is an $n$th power in $k_v$, thus $C/_{k_v} \cong ({C_d})/_{k_v}$, so $({C_d})/_{k_v}\neq \emptyset$. Finally, for $\mathfrak L_v = (\pi)$, by construction $k_v$ contains a root of $f(x)$, and so $({C_d})/_{k_v} \neq \emptyset$ for every twist of $C$. 
\par We take for $\ang \psi$ in Theorem \ref{TAHP} the unique subgroup of order $p$ in $Aut(C/\kb)$. Then $(C/\ang \psi) \cong C^{(N)}$. By Theorem \ref{TAHP}, only finitely many of the degree $p$-twists $C_d$, $d = \pi^p$ have $k$-rational points. \hfill $\square$
\subsection{Examples of the unconditional result} Examples satisfying the hypotheses of Cor \ref{unsup} include hyperelliptic curves of the form $y^2 = x^{2\ell}-A$, where $\ell \ge 5$ (or $\ell = 3,4$ and the genus 1 curve $y^2 = x^\ell-A$ has finitely many $k$-rational points) and $A$ is not an $2\ell$th power in $k$. In such cases, $y^2 = x^\ell - A$ has finitely many points, so by Thm. \ref{TAHP}, only finitely many quadratic twists of $y^2 = x^{2\ell}-A$ (with respect to the automorphism $\tau(x,y) = (x, \zeta_{\ell}y)$) have $k$-rational points.
\par  Additional examples include superelliptic curves defined by $y^4 = f(x)$ where $\deg(f) = 2m >4 $ (or $\deg(f) = 4$ and $y^2 = f(x)$ is of genus 1 and has only finitely many points $k$-rational points) and $f(x) \in k[x]$ has no $k$-rational roots. As the hyperelliptic curve $y^2 = f(x)$ has only finitely many rational points, there are only finitely many quadratic twists corresponding to the automorphism $(x,y^2) \mapsto (x,-y^2)$ which have $k$-rational points.

\end{document}